\documentclass[a4paper,12pt]{amsart}
\usepackage{amssymb}
\usepackage{graphicx}
\usepackage{ifthen}
\usepackage{xcolor}
\usepackage{cite}
\nonstopmode \numberwithin{equation}{section}
\newtheorem{thm}{Theorem}
\newtheorem{cor}{Corollary}
\newtheorem{lem}{Lemma}
\newtheorem{prop}{Proposition}

\newtheorem{conj}{Conjecture}

\theoremstyle{definition}
\newtheorem{defn}{Definition}[section]

\newtheorem{prob}[equation]{Problem}

\newenvironment{rem}{%
\bigskip
\noindent \textsl{{\bf Remark. }}}{\bigskip}
\newenvironment{rems}{%
\bigskip
\noindent \textsl{{\bf Remarks. }}}{\bigskip}

\newcounter {own}
\def\theown {\thesection       .\arabic{own}}

\newenvironment{pf}[1][]{%
 \vskip 3mm
 \noindent
 \ifthenelse{\equal{#1}{}}%
  {{\slshape Proof. }}%
  {{\slshape #1.} }%
 }%
{\qed\bigskip}

\newcounter{alphabet}

\newenvironment{Thm}[1][]{\refstepcounter{alphabet}%
\bigskip%
\noindent%
{\bf Theorem \Alph{alphabet}}%
\ifthenelse{\equal{#1}{}}{}{ (#1)}%
{\bf .} \itshape}{\vskip 8pt}

\newcommand{\IR}{{\mathbb R}}
\newcommand{\ID}{{\mathbb D}}
\newcommand{\IN}{{\mathbb N}}
\newcommand{\IC}{{\mathbb C}}

\newcommand{\IH}{{\mathcal H}}

\def\be{\begin{equation}}
\def\ee{\end{equation}}

\newcommand{\bee}{\begin{enumerate}}
\newcommand{\eee}{\end{enumerate}}

\newcommand{\blem}{\begin{lem}}
\newcommand{\elem}{\end{lem}}
\newcommand{\bthm}{\begin{thm}}
\newcommand{\ethm}{\end{thm}}
\newcommand{\bcor}{\begin{cor}}
\newcommand{\ecor}{\end{cor}}
\newcommand{\beg}{\begin{examp}}
\newcommand{\eeg}{\end{examp}}
\newcommand{\begs}{\begin{examples}}
\newcommand{\eegs}{\end{examples}}
\newcommand{\bdefe}{\begin{defn}}
\newcommand{\edefe}{\end{defn}}
\newcommand{\bprob}{\begin{prob}}
\newcommand{\eprob}{\end{prob}}
\newcommand{\bei}{\begin{itemize}}
\newcommand{\eei}{\end{itemize}}

\newcommand{\bcon}{\begin{conj}}
\newcommand{\econ}{\end{conj}}
\newcommand{\bcons}{\begin{conjs}}
\newcommand{\econs}{\end{conjs}}
\newcommand{\bprop}{\begin{prop}}
\newcommand{\eprop}{\end{prop}}
\newcommand{\br}{\begin{rem}}
\newcommand{\er}{\end{rem}}
\newcommand{\brs}{\begin{rems}}
\newcommand{\ers}{\end{rems}}
\newcommand{\bo}{\begin{obser}}
\newcommand{\eo}{\end{obser}}
\newcommand{\bos}{\begin{obsers}}
\newcommand{\eos}{\end{obsers}}
\newcommand{\bpf}{\begin{pf}}
\newcommand{\epf}{\end{pf}}
\newcommand{\ba}{\begin{array}}
\newcommand{\ea}{\end{array}}
\newcommand{\beq}{\begin{eqnarray}}
\newcommand{\beqq}{\begin{eqnarray*}}
\newcommand{\eeq}{\end{eqnarray}}
\newcommand{\eeqq}{\end{eqnarray*}}

\newcounter{minutes}
\divide\time by 60
\newcounter{hours}
\multiply\time by 60 \addtocounter{minutes}{-\time}
\begin{document}
\title{Estimates for generalized Bohr radii in one and higher dimensions}
\author{Nilanjan Das}
\address{Nilanjan Das, Theoretical Statistics and Mathematics Unit, Indian Statistical Institute Kolkata, Kolkata-700108, India.}
\email{nilanjand7@gmail.com}


\begin{abstract}
The generalized Bohr radius $R_{p, q}(X), p, q\in[1, \infty)$ for a complex Banach space $X$ was introduced by Blasco in 2010. In this article, we determine the exact value of $R_{p, q}(\IC)$ for the cases 
(i) $p, q\in[1, 2]$, (ii) $p\in (2, \infty), q\in [1, 2]$ and (iii) $p, q\in [2, \infty)$. Moreover, we consider an $n$-variable version $R_{p, q}^n(X)$ of the quantity $R_{p, q}(X)$ and determine (i) $R_{p, q}^n(\IH)$ for an infinite dimensional complex Hilbert space $\IH$, (ii) the precise asymptotic value of $R_{p, q}^n(X)$ as $n\to\infty$ for finite dimensional $X$. We also study the multidimensional analogue of a related concept called the $p$-Bohr radius, introduced by Djakov and Ramanujan in 2000. In particular, we obtain the asymptotic value of the $n$-dimensional $p$-Bohr radius for bounded complex-valued functions, and in the vector-valued case we provide a lower estimate for the same, which is independent of $n$. In a similar vein, we investigate in detail the multidimensional $p$-Bohr radius problem for functions with positive real part. Towards the end of this article, we pose one more generalization $R_{p, q}(Y, X)$ of $R_{p, q}(X)$-considering functions that map the open unit ball of another complex Banach space $Y$ inside the unit ball of $X$, and show that the existence of nonzero $R_{p, q}(Y, X)$ is governed by the geometry of $X$ alone.
\end{abstract}
\thanks{The author of this article is supported by a Research Associateship provided by the Stat-Math Unit of ISI Kolkata. .}

\maketitle
\pagestyle{myheadings}
\markboth{Estimates for generalized Bohr radii in one and higher dimensions}{N. Das}

\bigskip
\section{Introduction and the main results}\label{N2sec1}
The celebrated theorem of Harald Bohr \cite{Bohr} states (in sharp form) that for any holomorphic
self mapping $f(z)=\sum_{n=0}^\infty a_nz^n$ of the open unit disk $\ID$,
$$
\sum_{n=0}^\infty|a_n|r^n\leq 1
$$
for $|z|=r\leq 1/3$, and this quantity $1/3$ is the best possible. Inequalities of the above type are commonly known as \emph{Bohr inequalities} nowadays, and appearance of any such inequality in a result is generally termed as the occurrence of the \emph{Bohr phenomenon}.
This theorem was an outcome of Bohr's investigation on the \lq\lq absolute convergence problem\rq\rq of ordinary Dirichlet series of the form $\sum a_nn^{-s}$, and did not receive much attention until it was applied to answer a long-standing question in the realm of operator algebras in 1995 (cf. \cite{Dix}). Starting there, the Bohr phenomenon continues to be studied from several different aspects for the last two decades, for example: in certain abstract settings (cf. \cite{Aiz}), for ordinary and vector-valued Dirichlet series (see f.i. \cite{Bala, Def2}), for uniform algebras (see \cite{Paul2}), for free holomorphic functions (cf. \cite{Pop1}), for a Faber-Green condenser (see \cite{La}), for vector-valued functions (cf. \cite{Def, Ham, Ham1}), for Hardy space functions (see \cite{Bene}), and also for functions in several variables (see, for example \cite{BB2, Boas, Ga, Pop}). We also urge the reader to glance through the references of these above-mentioned articles to get a more complete picture of the recent developments in this area.
    
We will now concentrate on a variant of the Bohr inequality, introduced for the first time in \cite{Bla1} in order to investigate the Bohr phenomenon on Banach spaces. Let us start by defining an $n$-variable analogue of this modified inequality. For this purpose, we need to introduce some concepts. Let $\ID^n=\{(z_1, z_2, \cdots, z_n)\in\IC^n:\|z\|_\infty:=\max_{1\leq k\leq n}|z_k|<1\}$
be the open unit polydisk in the $n$-dimensional
complex plane $\IC^n$ and let $X$
be a complex Banach space. Any holomorphic function $f:\ID^n\to X$ can be expanded in the power series
\be\label{N2eq2}
f(z)=x_0+ \sum_{|\alpha|\in\IN} x_\alpha z^\alpha,\, x_\alpha\in X,
\ee
for $z\in\ID^n$.
Here and hereafter, we will use the standard multi-index notation: $\alpha$ denotes an $n$-tuple $(\alpha_1, \alpha_2,\cdots, \alpha_n)$ of
nonnegative integers, $|\alpha|:=\alpha_1+\alpha_2+\cdots+\alpha_n$, $z$ denotes an $n$-tuple $(z_1, z_2, \cdots, z_n)$ of
complex numbers, and $z^\alpha$ is the product $z_1^{\alpha_1}z_2^{\alpha_2}\cdots z_n^{\alpha_n}$. 
For $1\leq p,q<\infty$ and for any $f$ as in $(\ref{N2eq2})$ with $\|f\|_{H^\infty(\ID^n,X)}\leq 1$, we denote
$$
R_{p, q}^n(f, X)=\sup\left\{r\geq 0:\|x_0\|^p+\left(\sum_{k=1}^\infty\sum_{|\alpha|=k}\|x_\alpha z^\alpha\|\right)^q\leq 1\,\mbox{for all}\,z\in r\ID^n\right\},
$$
where $H^\infty(\ID^n,X)$ is the space of bounded holomorphic functions $f$ from $\ID^n$ to $X$ and $\|f\|_{H^\infty(\ID^n,X)}=\sup_{z\in\ID^n}\|f(z)\|$. We further define
$$
R_{p, q}^n(X)=\inf\left\{R_{p,q}^n(f, X):\|f\|_{H^\infty(\ID^n,X)}\leq 1\right\}.
$$
Following the notations of \cite{Bla1}, throughout this article we will use $R_{p, q}(f, X)$ for $R_{p, q}^1(f, X)$ and $R_{p, q}(X)$ for $R_{p, q}^1(X)$.
Clearly, $R_{1, 1}(\IC)=1/3$. The reason of reshaping the original Bohr inequality in the above fashion becomes clear from \cite[Theorem 1.2]{Bla1}, which shows that the notion of the classical Bohr phenomenon is not very useful for $\mbox{dim}(X)\geq 2$. For a given pair of $p$ and $q$ in $[1, \infty)$, it is known from the results of \cite{Bla1} that depending on $X$, $R_{p, q}(X)$ may or may not be zero. A characterization theorem in this regard has further been established in \cite{BB}. However, the question of determination of the exact value of $R_{p, q}(X)$ is challenging, and to the best of our knowledge there is lack of progress on this problem-even for $X=\IC$. In fact, only known optimal result in this direction is the following:
\be\label{N2eq3}
R_{p, 1}(\IC)=\frac{p}{2+p}
\ee
for $1\leq p\leq 2$ (cf. \cite[Proposition 1.4]{Bla1}), along with rather recent generalizations of $(\ref{N2eq3})$ (see, for example \cite{Liu}). This motivates us to address this problem in the first theorem of this article.
\bthm\label{N2thm1}
Given $p, q\in[1, \infty)$, let us denote
$$
A_{p, q}(a)=
\frac{(1-a^p)^{\frac{1}{q}}}{1-a^2+a(1-a^p)^{\frac{1}{q}}},\, a\in [0, 1)
$$
and 
$$
S_{p, q}(a)=\left(\frac{(1-a^p)^{\frac{2}{q}}}{1-a^2+(1-a^p)^{\frac{2}{q}}}\right)^{\frac{1}{2}},\, a\in [0, 1).
$$
Also let $\widehat{a}$ be the unique root in $(0, 1)$ of the equation 
\be\label{N2eq4}
x^p+x^q=1.
\ee
Then
$$
R_{p, q}(\IC)=
\begin{cases}
\inf\limits_{a\in[\widehat{a}, 1)}A_{p, q}(a)\, \,\emph{if}\,\, p, q\in[1,2],\\
\min\left\{(1/\sqrt{2}), \inf\limits_{a\in[\widehat{a}, 1)}A_{p, q}(a)\right\}\,\, \emph{if}\,\, p\in(2, \infty)\, \emph{and}\, q\in [1, 2],\\
1/\sqrt{2}\,\,\emph{if}\,\, p, q\in [2, \infty).
\end{cases}
$$
For $p\in [1, 2]$ and $q\in (2, \infty)$,
$R_{2, q}(\IC)=1/\sqrt{2}$,
$R_{p, q}(\IC)=\inf_{a\in[\widehat{a}, 1)}A_{p, q}(a)$ if $p<2$ and in addition the inequality
\be\label{N2eq5}
q\widehat{a}^2+p\widehat{a}^{p+2}\leq p\widehat{a}^p+q\widehat{a}^{p+2}
\ee
is satisfied. In all other scenarios, we have, in general
\be\label{N2eq15}
0<\inf_{a\in[0, 1)}S_{p, q}(a)\leq R_{p, q}(\IC)\leq \frac{1}{\sqrt{2}}.
\ee
\ethm
\brs
\textbf{(a)}~A closer look at the proof of Theorem \ref{N2thm1} reveals that the conclusions of this theorem remain unchanged if the interval $[1, 2]$ is replaced by $(0, 2]$ everywhere in its statement. However, doing so includes cases where positive Bohr radius is nonexistent-for example, $R_{p, q}(\IC)=\inf_{a\in[\widehat{a}, 1)}A_{p, q}(a)\leq\lim_{a\to 1-}A_{p, q}(a)=0$ if $0<q<1$. Therefore, throughout this paper we stick to the assumption $p, q\geq 1$.

\noindent
\textbf{(b)}~Following methods similar to the proof of Theorem \ref{N2thm1}, it is easy to see that
for any given complex Hilbert space $\IH$ with dimension at least $2$, the following statements are true:

\noindent
(i) For $p, q\in[2, \infty)$, $R_{p, q}(\IH)=1/\sqrt{2}$.

\noindent
(ii) For $p\in[1, 2)$ and $q\in [2, \infty)$, the inequalities $(\ref{N2eq15})$ are satisfied with $R_{p, q}(\IC)$ replaced by $R_{p, q}(\IH)$.

\noindent
Note that the assumption $q\geq 2$ is justified by \cite[Corollary 4]{BB}.
Later in Theorem \ref{N2thm3}, we obtain a more complete result for $\mbox{dim}(\IH)=\infty$. 
\ers

We now turn our attention to the Bohr radius $R_{p, q}^n(X)$, where $X$ is a complex Banach space. The first question we encounter is the identification of the Banach spaces $X$ with $R_{p, q}^n(X)>0$; which is in fact equivalent to the one-dimensional version of the same problem.
\bprop\label{N2thm2}
For any given $n\in\IN$ and $p, q\in[1, \infty)$, $R_{p, q}^n(X)>0$ for some complex Banach space $X$ if and only if $R_{p, q}(X)>0$ for the same Banach space $X$.
\eprop
Note that from \cite[Theorem 1]{BB} it is known that $R_{p, q}(X)>0$ if and only if there exists a constant $C$ such that
\be\label{N2eq6}
\Omega_X(\delta)\leq C\left((1+\delta)^q-(1+\delta)^{q-p}\right)^{1/q} 
\ee
for all $\delta\geq 0$. We mention here that for any $\delta\geq 0$, $\Omega_X(\delta)$ is defined to be the supremum of $\|y\|$ taken over all $x, y\in X$ such that $\|x\|=1$ and $\|x+zy\|\leq 1+\delta$ for all $z\in\ID$ (see \cite{Glo}).
Now, in view of the above discussion, it looks appropriate to consider the Bohr phenomenon, i.e. studying $R_{p, q}^n(X)$ for particular Banach spaces $X$. We resolve this problem completely for $X=\IH$-a complex Hilbert space of infinite dimension. While this question remains open for $\mbox{dim}(\IH)<\infty$, we succeed in determining the correct asymptotic behaviour of $R_{p, q}^n(X)$ as $n\to\infty$ for any finite dimensional complex Banach space $X$ with $R_{p, q}(X)>0$.
\bthm\label{N2thm3}
For any given $n\in\IN$, $p\in[1, \infty)$, $q\in[2, \infty)$ and for any infinite dimensional complex Hilbert space $\IH$:
$$
R_{p, q}^n(\IH)=\inf_{a\in[0,1)}\left(1-(1-(S_{p,q}(a))^2)^{\frac{1}{n}}\right)^{\frac{1}{2}},
$$
$S_{p, q}(a)$ as defined in the statement of Theorem \ref{N2thm1}. For any complex Banach space $X$ with $\emph{dim}(X)<\infty$ and with $R_{p, q}(X)>0$, we have $$
\lim_{n\to\infty}R_{p, q}^n(X)\sqrt{\frac{n}{\log n}}=1.
$$
\ethm
At this point, we like to discuss another interesting related concept called the $p$-Bohr radius. First we pose an $n$-variable version of the definition of $p$-Bohr radius given in \cite{Bla2}. For any $p\in[1, \infty)$ and for any complex Banach space $X$, we denote
$$
r_p^n(f, X)=\sup\left\{r\geq 0:\|x_0\|^p+\sum_{k=1}^\infty\sum_{|\alpha|=k}\|x_\alpha z^\alpha\|^p\leq 1\,\mbox{for all}\,z\in r\ID^n\right\},
$$
where $f$ is as given in $(\ref{N2eq2})$ with $\|f\|_{H^\infty(\ID^n, X)}\leq 1$, and then define the $n$-dimensional $p$-Bohr radius of $X$ by
$$
r_p^n(X)=\inf\left\{r_p^n(f, X):\|f\|_{H^\infty(\ID^n,X)}\leq 1\right\}.
$$
Again, following the notations of \cite{Bla2}, we will write $r_p(f, X)$ for $r_p^1(f, X)$ and $r_p(X)$ for $r_p^1(X)$. Clearly, for $X=\IC$ one only needs to consider $p\in[1,2)$, as $r_p^n(\IC)=1$ for all $p\geq 2$ and for any $n\in\IN$. The quantities $r_p(\IC)$ and $r_p^n(\IC)$ were first considered in \cite{Dja}. Unlike $R_{p, q}(\IC)$, precise value of $r_p(\IC)$ has already been obtained in \cite{Kay1}. We make further progress by determining the asymptotic behaviour of $r_p^n(\IC)$ for all $p\in (1,2)$ (the cases $p=1$ is already resolved) in the first half of Theorem \ref{N2thm4}.

On the other hand, to get a nonzero value of $r_p^n(X)$ where $\mbox{dim}(X)\geq 2$, one necessarily has to consider $p\geq 2$ and work  with $p$-uniformly $PL$-convex complex Banach spaces $X$. A complex Banach space $X$ is said to be $p$-uniformly $PL$-convex ($2\leq p<\infty$) if there exists a constant $\lambda>0$ such that
\be\label{N2eq7}
\|x\|^p+\lambda\|y\|^p\leq\frac{1}{2\pi}\int_0^{2\pi}\|x+e^{i\theta}y\|^pd\theta
\ee
for all $x, y\in X$. Denote by $I_p(X)$ the supremum of all $\lambda$ satisfying $(\ref{N2eq7})$. Now if we assume $r_p^n(X)>0$ for some $n\in\IN$ then evidently $r_p(X)>0$ (as any member of $H^\infty(\ID, X)$ can be considered as a member of $H^\infty(\ID^n, X)$ as well), and therefore \cite[Theorem 1.10]{Bla2} asserts that $X$ is $p$-uniformly $\IC$-convex; which is equivalent to saying $X$ is $p$-uniformly $PL$-convex.
The second half of our upcoming theorem shows that for any $p$-uniformly $PL$-convex complex Banach space $X$ ($p\geq 2$) with $\mbox{dim}(X)\geq 2$, the Bohr radius $r_p^n(X)>0$ for all $n\in\IN$ and unlike $r_p^n(\IC)$ or $R_{p, q}^n(X)$, $r_p^n(X)$ does not converge to $0$ as $n\to\infty$.
\bthm\label{N2thm4}
For any $p\in (1,2)$ and $n>1$, we have
$$
r_p^n(\IC)
\thicksim\left(\frac{\log n}{n}\right)^{\frac{2-p}{2p}}.
$$
For any $p$-uniformly $PL$-convex $(p\geq 2)$ complex Banach space $X$ with $\mbox{dim}(X)\geq 2$, we have
$$
\left(\frac{I_p(X)}{2^p+I_p(X)}\right)^{\frac{2}{p}}\leq r_p^n(X)\leq 1
$$
for all $n\in\IN$.
\ethm
We clarify that for any two sequences $\{p_n\}$ and $\{q_n\}$ of positive real numbers, we write $p_n\thicksim q_n$ if there exist constants $C, D>0$ such that
$Cq_n\leq p_n\leq Dq_n$ for all $n>1$.

Let us now focus on the Bohr phenomenon for complex-valued functions that are not necessarily bounded-in particular, for functions that map inside the right half-plane. The first result in this direction was proved by Aizenberg et. al. \cite{Aiz}, which states that for any holomorphic function $f(z)=\sum_{n=0}^\infty a_nz^n$ defined on $\ID$ that satisfies $\mbox{Re}(f(z))>0$ for $z\in\ID$, $f(0)>0$, the inequality 
$$
\sum_{n=0}^\infty|a_n|r^n\leq 2f(0)
$$
holds for $|z|=r\leq 1/3$, and this constant $1/3$ is the best possible. It is clear that without loss of generality, $f(0)=1$ could be assumed in the above inequality. Extensions of this result in several variable framework can be found in \cite{Aiz2}, and in \cite{Dja} $p$-Bohr radius for such functions was considered in single variable setting. Our aim is to blend these two approaches together. To this end, we introduce the following quantities. For $p>0$, we denote
$$
H_p^n(f)=\sup\left\{r\geq 0:\frac{1}{2}\left(|x_0|^p+\sum_{k=1}^\infty\sum_{|\alpha|=k}|x_\alpha z^\alpha|^p\right)^{\frac{1}{p}}\leq 1\,\mbox{for all}\,z\in r\ID^n\right\},
$$
where $f:\ID^n\to\IC$ is as given in $(\ref{N2eq2})$, with complex coefficients $x_\alpha$ and with $\mbox{Re}(f(z))>0$ for all $z\in\ID^n$, $f(0)=1$. We also define 
$$
H_p^n=\inf\left\{H_p^n(f):\mbox{Re}(f(z))>0\,\,\mbox{for all}\,\,z\in\ID^n, f(0)=1\right\}.
$$
Note that $p$-Bohr radius of the above type makes sense for $p\in(0,1)$, which is not the case for $r_p^n(\IC)$. Clearly, $H_1^1=1/3$ and the value of $H_p^1$ is known from \cite{Dja} for all other $p>0$. As a step forward, in the following theorem we describe the behaviour of $H_p^n$ in different ranges of $p$.
\bthm\label{N2thm5}
For any $n>1$:
$$
H_p^n=\left(\frac{2^p-1}{2^{p+1}-1}\right)^{\frac{1}{p}}
$$
for $p\in [2, \infty)$, and 
$$
H_p^n
\thicksim\left(\frac{\log n}{n}\right)^{\frac{2-p}{2p}}
$$
for $p\in(0, 2)$.
\ethm

In Section \ref{N2sec2}, we will give the proofs of all the results stated so far. Finally, it is natural to think of extending the concept of $R_{p,q}(X)$ by considering functions from a domain in another complex Banach space $Y\neq\IC$ to the complex Banach space $X$. In Section \ref{N2sec3}, we consider one such generalization of $R_{p, q}(X)$ and observe that the occurrence of the Bohr phenomenon remains independent of the choice of $Y$. 
\section{Proofs of the main results}\label{N2sec2}
We start by recalling the following result of Bombieri (cf. \cite{Bo}), which is at the heart of the proof of our Theorem \ref{N2thm1}.

\begin{Thm}\label{N2TheoA}
For any holomorphic self mapping $f(z)=\sum_{n=0}^\infty a_nz^n$ of the open unit disk $\ID$:
$$
\sum_{n=1}^\infty |a_n|r^n\leq
\begin{cases}
\frac{r(1-a^2)}{1-ar}\, \,\emph{for}\,\, r\leq a,\\
\frac{r\sqrt{1-a^2}}{\sqrt{1-r^2}}\,\, \emph{for}\,\,r\in [0,1)\,\,\emph{in general},
\end{cases}
$$
where $|z|=r$ and $|a_0|=a$.
\end{Thm}
It should be mentioned that the above result is not recorded in the present form in \cite{Bo}. For a direct derivation of the first inequality in Theorem A, see the proof of Theorem 9 from \cite{BB1}. The second inequality is an easy consequence of the Cauchy-Schwarz inequality combined with the fact $\sum_{n=1}^\infty|a_n|^2\leq 1-|a_0|^2$.
\bpf[\bf Proof of Theorem \ref{N2thm1}]
Given a holomorphic function $f(z)=\sum_{n=0}^\infty a_nz^n$ mapping $\ID$ inside $\ID$, a straightforward application of Theorem A yields
\be\label{N2eq8}
|a_0|^p+\left(\sum_{n=1}^\infty |a_n|r^n\right)^q\leq
\begin{cases}
a^p+(1-a^2)^q\left(\frac{r}{1-ar}\right)^q\, \,\textrm{for}\,\, r\leq a,\\
a^p+(1-a^2)^{\frac{q}{2}}\left(\frac{r}{\sqrt{1-r^2}}\right)^q\,\, \textrm{for}\,\,r\in [0,1).
\end{cases}
\ee
Now, 
$$
a^p+(1-a^2)^q\left(\frac{r}{1-ar}\right)^q\leq 1
$$
whenever $r\leq A_{p, q}(a)$. A little calculation reveals that $A_{p,q}(a)\leq a$ whenever $a^p+a^q\geq 1$, i.e. whenever $a\geq\widehat{a}$, $\widehat{a}$ being the root of the equation $(\ref{N2eq4})$. Thus from $(\ref{N2eq8})$, it is clear that 
\be\label{N2eq22}
|a_0|^p+\left(\sum_{n=1}^\infty |a_n|r^n\right)^q\leq 1
\ee
for $r\leq\inf_{a\in[\widehat{a}, 1)}A_{p, q}(a)$, provided that $a\geq \widehat{a}$. On the other hand, 
$$
a^p+(1-a^2)^{\frac{q}{2}}\left(\frac{r}{\sqrt{1-r^2}}\right)^q\leq 1
$$
for $r\leq S_{p, q}(a)$, i.e. the inequality $(\ref{N2eq22})$ remains valid for $r\leq\inf_{a\in[0, \widehat{a}]}S_{p, q}(a)$, provided that $a\leq \widehat{a}$. Therefore, we conclude that for any given $p, q\in[1, \infty)$
\be\label{N2eq9}
R_{p, q}(\IC)\geq \min\left\{\inf_{a\in[0, \widehat{a}]}S_{p, q}(a), \inf_{a\in[\widehat{a}, 1)}A_{p, q}(a)\right\}.
\ee
We also record some other facts which we will need to use later.
Observe that for all $p, q\in[1, \infty)$
$$
S_{p, q}(a)=\sqrt{\frac{T(a)}{1+T(a)}}\,\, \mbox{where}\,\,T(a)=\frac{(1-a^p)^{\frac{2}{q}}}{1-a^2},
$$
and therefore
$$
S_{p, q}^\prime(a)
=\frac{T^\prime(a)}{2\sqrt{T(a)(1+T(a))^3}}
$$
for $a\in(0,1)$, where
\be\label{N2eq10}
T^\prime(a)
=\frac{2a^{p-1}T(a)}{1-a^p}\left(\frac{a^2(1-a^p)}{a^p(1-a^2)}-\frac{p}{q}\right).
\ee
Setting $y=1/a$ for convenience, we write
$$
\frac{a^2(1-a^p)}{a^p(1-a^2)}
=\frac{y^p-1}{y^2-1}=P(y)
$$
defined on $(1,\infty)$. Note that
\be\label{N2eq11}
\frac{d}{da}P(y)=P^\prime(y)\frac{dy}{da}
=-y^3\frac{py^p-py^{p-2}-2y^p+2}{(y^2-1)^2},
\ee
and that
\be\label{N2eq12}
Q^\prime(y)=y^{p-3}(y^2-1)p(p-2),
\ee
where $Q(y)=py^p-py^{p-2}-2y^p+2$.

Further, observe that
for the disk automorphisms $\phi_a(z)=(a-z)/(1-az)$, $z\in\ID, a\in[\widehat{a}, 1)$, $R_{p, q}(\phi_a, \IC)=A_{p, q}(a)$, and hence $R_{p, q}(\IC)\leq\inf_{a\in[\widehat{a}, 1)}A_{p, q}(a)$. Also, for $\xi(z)=z\phi_{1/\sqrt{2}}(z)$, $z\in\ID$ we have
$R_{p, q}(\xi, \IC)=1/\sqrt{2}$. Combining these two facts, we write
\be\label{N2eq14}
R_{p, q}(\IC)\leq\min\left\{(1/\sqrt{2}), \inf_{a\in[\widehat{a}, 1)}A_{p, q}(a)\right\}.
\ee

\vspace{2pt}
\noindent We now deal with the problem case by case.

\vspace{3pt}
\noindent\underline{Case $p, q\in [1,2]$}:
Let us start with $p<2$.
From $(\ref{N2eq12})$, it is evident that 
$Q^\prime(y)<0$ for $p<2$, and hence
$Q(y)<Q(1)=0$ for all $y\in(1, \infty)$. Thus from $(\ref{N2eq11})$ it is clear that $P(y)$ is strictly increasing in $(0,1)$ with respect to $a$. Consequently,
for all $y\in(1, \infty)$
\be\label{N2eq13}
P(y)<\lim_{a\to 1-}P(y)=\frac{p}{2},
\ee
and using the above estimate in $(\ref{N2eq10})$ gives, for all $a\in(0,1)$
$$
T^\prime(a)
<\frac{2a^{p-1}T(a)}{1-a^p}\left(\frac{p}{2}-\frac{p}{q}\right)
\leq 0,
$$
as $q\leq 2$. Therefore $S_{p, q}(a)$ is strictly decreasing in $(0,1)$, and after some calculations we have, as a consequence:
$$
\inf_{a\in[0, \widehat{a}]}S_{p, q}(a)=S_{p, q}(\widehat{a})=A_{p, q}(\widehat{a})\geq  \inf_{a\in[\widehat{a}, 1)}A_{p, q}(a).
$$
Hence from $(\ref{N2eq9})$, we have $R_{p,q}(\IC)\geq\inf_{a\in[\widehat{a}, 1)}A_{p, q}(a)$. For $p=2$, if $q<2$ then $T^\prime(a)<0$ for all $a\in(0,1)$, which (as in the case $p<2$) again gives $R_{2,q}(\IC)\geq\inf_{a\in[\widehat{a}, 1)}A_{2, q}(a)$. Otherwise, if $p=q=2$ then $\widehat{a}=1/\sqrt{2}$, and for all $a\in[0, 1)$ we get 
$$
S_{2, 2}(a)=1/\sqrt{2}=\inf_{a\in[\widehat{a}, 1)}A_{2, 2}(a).
$$
Therefore, for all $p, q\in[1,2]$ we have
$R_{p,q}(\IC)\geq\inf_{a\in[\widehat{a}, 1)}A_{p, q}(a)$,
and from $(\ref{N2eq14})$ it is known that $R_{p, q}(\IC)\leq\inf_{a\in[\widehat{a}, 1)}A_{p, q}(a)$. This completes the proof for this case.

\vspace{3pt}
\noindent\underline{Case $p\in(2, \infty), q\in [1,2]$}:
From $(\ref{N2eq12})$ it is clear that $Q^\prime(y)>0$ for $p>2$, and therefore $Q(y)>Q(1)=0$ for all $y\in(1, \infty)$. It follows from $(\ref{N2eq11})$ that $P(y)$ is strictly decreasing in $(0,1)$ with respect to $a$. Thus, for $q<2$, the value of the quantity
$$
P(y)-\frac{p}{q}=\frac{a^2(1-a^p)}{a^p(1-a^2)}-\frac{p}{q}
$$
decreases from 
$$
\lim_{a\to 0+}(P(y)-(p/q))=+\infty\,\, \mbox{to}\,\, \lim_{a\to 1-}(P(y)-(p/q))=p((1/2)-(1/q))<0,
$$
i.e. $P(y)-(p/q)>0$ in $(0, b_1)$ and $P(y)-(p/q)<0$ in $(b_1, 1)$ for some $b_1\in (0,1)$, where $P(b_1)=(p/q)$.
As a consequence, $T^\prime(a)=0$ only for $a=0, b_1$, and $T^\prime(a)>0$ in $(0, b_1)$, $T^\prime(a)<0$ in $(b_1, 1)$. Hence, $S_{p, q}(a)$ strictly increases in $(0, b_1)$, and then strictly decreases in $(b_1, 1)$, which implies
$$
\inf_{a\in[0, \widehat{a}]}
S_{p, q}(a)=\min\left\{S_{p, q}(0), S_{p, q}(\widehat{a})\right\}
=\min\left\{(1/\sqrt{2}), A_{p, q}(\widehat{a})\right\}.
$$
Moreover, from the proof of the case $p, q\in[2, \infty)$, we have $R_{p, 2}(\IC)=1/\sqrt{2}$. These two facts combined with $(\ref{N2eq9})$ readily yields
$$
R_{p, q}(\IC)\geq
\min\left\{(1/\sqrt{2}), \inf\limits_{a\in[\widehat{a}, 1)}A_{p, q}(a)\right\},
$$
and making use of $(\ref{N2eq14})$ we arrive at our desired conclusion.

\vspace{3pt}
\noindent\underline{Case $p, q\in [2, \infty)$}: Applying $(\ref{N2eq14})$ of this paper, $(1.9)$ from \cite{Bla1} and \cite[Remark 1.2]{Bla2} together, the proof follows immediately from the observation:
$$
(1/\sqrt{2})\geq R_{p, q}(\IC)\geq R_{2, 2}(\IC)\geq (1/\sqrt{2})r_2(\IC)=1/\sqrt{2}.
$$

\vspace{3pt}
\noindent\underline{Case $p\in [1,2], q\in (2, \infty)$}: The fact that $R_{2, q}(\IC)=1/\sqrt{2}$ is evident from the proof of the case $p, q\in[2, \infty)$. Further, as we have already seen, from $(\ref{N2eq8})$ it is clear that 
the inequality $(\ref{N2eq22})$ holds
for $r\leq S_{p, q}(a), a\in [0,1)$, and therefore for $r\leq\inf_{a\in[0,1)}S_{p ,q}(a)$. From this and $(\ref{N2eq14})$, we have $(\ref{N2eq15})$ as an immediate consequence. The assertion $\inf_{a\in[0,1)}S_{p, q}(a)>0$ is validated from the fact that $S_{p, q}(a)\neq 0$ for all $a\in[0,1)$ and that $\lim_{a\to 1-}S_{p, q}(a)=1$. Now we will show that the imposition of the additional condition $(\ref{N2eq5})$ gives an optimal value for $R_{p, q}(\IC)$. We know that for $p<2$, $P(y)$ is strictly increasing in $(0,1)$ with respect to $a$, and as a result $P(y)-(p/q)$ increases from 
$$
\lim_{a\to 0+}(P(y)-(p/q))=-p/q\,\, \mbox{to}\,\, \lim_{a\to 1-}(P(y)-(p/q))=p((1/2)-(1/q))>0,
$$
i.e. $P(y)-(p/q)<0$ in $(0, b_2)$ and $P(y)-(p/q)>0$ in $(b_2, 1)$ for some $b_2\in (0,1)$, where $P(b_2)=(p/q)$.
As a consequence, $T^\prime(a)=0$ only for $a=0, b_2$, and $T^\prime(a)<0$ in $(0, b_2)$, $T^\prime(a)>0$ in $(b_2, 1)$. Hence, $S_{p, q}(a)$ strictly decreases in $(0, b_2)$, and then strictly increases in $(b_2, 1)$. Now if we assume the condition $(\ref{N2eq5})$ in addition, it is equivalent to saying that $T^\prime(\widehat{a})\leq 0$, i.e. $\widehat{a}\leq b_2$. Thus, $\inf_{a\in[0, \widehat{a}]}S_{p, q}(a)=S_{p, q}(\widehat{a})=A_{p, q}(\widehat{a})$. Consequently, from $(\ref{N2eq9})$ we get
$R_{p, q}(\IC)\geq\inf_{a\in[\widehat{a}, 1)}A_{p, q}(a)$, which completes our proof for this case.
\epf

\bpf[\bf Proof of Proposition \ref{N2thm2}]
As any holomorphic function $f:\ID\to X$ can also be considered as a holomorphic function from $\ID^n$ to $X$, it immediately follows that $R_{p, q}^n(X)>0$ for any $n\in\IN$ implies $R_{p, q}(X)>0$. Thus we only need to establish the converse. Any holomorphic $f:\ID^n\to X$ with an expansion $(\ref{N2eq2})$ can be written as \be\label{N2eq16}
f(z)=x_0+\sum_{k=1}^\infty P_k(z), z\in\ID^n
\ee
where $P_k(z):=\sum_{|\alpha|=k}x_\alpha z^\alpha$. Thus for any fixed $z_0\in\mathbb{T}^n$-the $n$-dimensional torus, we have
\be\label{N2eq23}
g(u):=f(uz_0)=x_0+\sum_{k=1}^\infty P_k(z_0)u^k:\ID\to X
\ee
is holomorphic, and if $\|f\|_{H^\infty(\ID^n, X)}\leq 1$ then $\|g\|_{H^\infty(\ID, X)}\leq 1$. Hence, starting with the assumption $R_{p, q}(X)=R>0$, we have
$\|P_k(z_0)\|\leq(1/R^{k})(1-\|x_0\|^p)^{1/q}$,
and since $z_0$ is arbitrary, we conclude that 
$\sup_{z\in\mathbb{T}^n}\|P_k(z)\|
\leq(1/R^{k})(1-\|x_0\|^p)^{1/q}$
for any $k\in\IN$. Therefore, for a given $k\in\IN$ and for any $\alpha$ with $|\alpha|=k$, we have
\begin{align*}
\|x_\alpha\|&=\left\|\frac{1}{(2\pi i)^n}\int_{|z_1|=1}\int_{|z_2|=1}\cdots\int_{|z_n|=1}\frac{P_k(z)}{z^{\alpha+1}}dz_ndz_{n-1}\cdots dz_1\right\|\\
&\leq\sup_{z\in\mathbb{T}^n}\|P_k(z)\|\leq\frac{1}{R^{k}}(1-\|x_0\|^p)^{\frac{1}{q}}.
\end{align*}
As a result we have, for all $r<R$:
$$
\|x_0\|^p+\left(\sum_{k=1}^\infty r^k\sum_{|\alpha|=k}\|x_\alpha\|\right)^q
\leq \|x_0\|^p+(1-\|x_0\|^p)\left(\left(\frac{R}{R-r}\right)^n-1\right)^q,
$$
which is less than or equal to 1 whenever 
$r\leq R\left(1-\left(1/2\right)^{1/n}\right)$,
thereby asserting $R_{p, q}^n(X)>0$.
\epf
\bpf[\bf Proof of Theorem \ref{N2thm3}]
(i) Before we start proving the first part of this theorem, note that the choice of $q\in[2, \infty)$ is again justified due to Proposition \ref{N2thm2} and \cite[Corollary 4]{BB}. Now, given a holomorphic $f:\ID^n\to\IH$ with an expansion $(\ref{N2eq2})$ and with $\|f(z)\|\leq 1$ for all $z\in\ID^n$, 
we have 
$\|x_0\|^2+\sum_{k=1}^\infty\sum_{|\alpha|=k}\|x_\alpha\|^2\leq 1$.
Taking $z\in r\ID^n$ and using this inequality, we obtain
\begin{align*}
 \|x_0\|^p+\left(\sum_{k=1}^\infty \sum_{|\alpha|=k}\|x_\alpha z^\alpha\|\right)^q
 &\leq\|x_0\|^p+\left(\sum_{k=1}^\infty \sum_{|\alpha|=k}\|x_\alpha\|^2\right)^{\frac{q}{2}}\left(\sum_{k=1}^\infty \sum_{|\alpha|=k}|z^\alpha|^2\right)^{\frac{q}{2}}\\
 &\leq\|x_0\|^p+(1-\|x_0\|^2)^{\frac{q}{2}}\left(\sum_{k=1}^\infty {n+k-1 \choose k} r^{2k}\right)^{\frac{q}{2}}\\
 &=\|x_0\|^p+(1-\|x_0\|^2)^{\frac{q}{2}}\left(\frac{1}{(1-r^2)^n}-1\right)^{\frac{q}{2}},
\end{align*}
which is less than or equal to 1 if 
\be\label{N2eq17}
r\leq\left(1-(1-(S_{p,q}(\|x_0\|))^2)^{\frac{1}{n}}\right)^{\frac{1}{2}},
\ee
and therefore,
\be\label{N2eq27}
R_{p, q}^n(\IH)\geq \inf_{a\in[0,1)}\left(1-(1-(S_{p,q}(a))^2)^{\frac{1}{n}}\right)^{\frac{1}{2}}.
\ee
As the quantity on the right hand side of the inequality $(\ref{N2eq17})$ becomes $\sqrt{1-(1/2)^{1/n}}$ at $x_0=0$ and converges to $1$ as $\|x_0\|\to 1-$, we 
conclude that the infimum in the inequality $(\ref{N2eq27})$ is attained at some $b_3\in[0, 1)$.
Since every Hilbert space $\IH$ has an orthonormal basis and in our case, $\mbox{dim}(\IH)=\infty$, we can choose a countably infinite set $\{e_\alpha\}_{|\alpha|\in\IN\cup\{0\}}$ of orthonormal vectors in $\IH$.
Setting 
$r_3=(1-(1-(S_{p,q}(b_3))^2)^{\frac{1}{n}})^{\frac{1}{2}}$,
we construct
$$
\chi(z):=b_3e_0+\frac{1-b_3^2}{(1-b_3^p)^{\frac{1}{q}}}\sum_{k=1}^\infty r_3^k\left(\sum_{|\alpha|=k}z^\alpha e_\alpha\right):\ID^n\to\IH,
$$
which satisfies $\|\chi(z)\|\leq 1$ for all $z\in\ID^n$, and $r_3=R_{p, q}^n(\chi, \IH)\geq R_{p, q}^n(\IH)$. This completes the proof for the first part of this theorem.

\vspace{3pt}
\noindent
(ii) Since $R_{p, q}(X)>0$, from the newly defined function $g(u)$ in $(\ref{N2eq23})$ we have
\be\label{N2eq18}
\sup_{z\in\mathbb{T}^n}\|P_k(z)\|\leq 2C\left((2-\|x_0\|)^q-(2-\|x_0\|)^{q-p}\right)^{1/q}
\ee
for any $k\in\IN$, $C$ being the constant for which
$(\ref{N2eq6})$ is satisfied (see inequalities $(2.2)$ and $(2.3)$ from \cite{BB}). Now as $X$ is finite dimensional, it is known that there exists another constant $D$ such that
$$
\left(\sum_{|\alpha|=k}\|x_\alpha\|^{\frac{2k}{k+1}}\right)^{\frac{k+1}{2k}}\leq D\sup_{\phi\in B_{X^*}}\left(\sum_{|\alpha|=k}|\phi(x_\alpha)|^{\frac{2k}{k+1}}\right)^{\frac{k+1}{2k}}
$$
for all $k\in\IN$, $B_{X^*}$ being the open unit ball in the dual space $X^*$ (see Proposition 2.3 and Theorem 2.8 from \cite{Die}). From \cite[Theorem 1.1]{Bay}, we further get that for any $\epsilon >0$, there exists $\mu>0$ such that
$$
\sup_{\phi\in B_{X^*}}\left(\sum_{|\alpha|=k}|\phi(x_\alpha)|^{\frac{2k}{k+1}}\right)^{\frac{k+1}{2k}}
\leq\mu(1+\epsilon)^k\sup_{\phi\in B_{X^*}}\sup_{z\in\ID^n}\left|\phi\left(P_k(z)\right)\right|
=\mu(1+\epsilon)^k\sup_{z\in\mathbb{T}^n}\|P_k(z)\|
$$
for all $k\geq 1$. Combining the above two inequalities and the inequality $(\ref{N2eq18})$ appropriately, it follows that
\begin{align*}
\left(\sum_{k=1}^\infty r^k\sum_{|\alpha|=k}\|x_\alpha\|\right)^q
&\leq\left(\sum_{k=1}^\infty r^k\left(\sum_{|\alpha|=k}\|x_\alpha\|^{\frac{2k}{k+1}}\right)^{\frac{k+1}{2k}}
{n+k-1 \choose k}^{\frac{k-1}{2k}}\right)^q\\
&\leq X\left(\sum_{k=1}^\infty r^k(1+\epsilon)^k{n+k-1 \choose k}^{\frac{k-1}{2k}}\right)^q
\end{align*}
where
$X=\mu^q C_1^q\left((2-\|x_0\|)^q-(2-\|x_0\|)^{q-p}\right)$,
$C_1=2CD$. Hence for $z\in r\ID^n$, the inequality
$$
\|x_0\|^p+\left(\sum_{k=1}^\infty\sum_{|\alpha|=k}\|x_\alpha z^\alpha\|\right)^q\leq 1
$$
is satisfied if 
\be\label{N2eq19}
\left(\frac{X}{1-\|x_0\|^p}\right)^{\frac{1}{q}}\left(\sum_{k=1}^\infty r^k(1+\epsilon)^k{n+k-1 \choose k}^{\frac{k-1}{2k}}\right)\leq 1. 
\ee
Now, analyzing the function $f_1(t)=((2-t)^p-1)/(1-t^p), t\in[0,1)$ (see \cite[pp. 554-555] {BB}) we see that $f_1(t)\leq 2^p-1$ for all $t\in[0,1)$, and hence
$$
\frac{X}{1-\|x_0\|^p}
\leq
\begin{cases}
\mu^qC_1^q 2^{q-p}(2^p-1)\,\,\mbox{if}\,\,q\geq p,\\
\mu^qC_1^q(2^p-1)\,\,\mbox{if}\,\,q\leq p.
\end{cases}
$$
Thus, inequality $(\ref{N2eq19})$ is satisfied if 
$$
C_2\left(\sum_{k=1}^\infty r^k(1+\epsilon)^k{n+k-1 \choose k}^{\frac{k-1}{2k}}\right)\leq 1,
$$
where $C_2$ is a new constant depending on $\mu, p, q$ and the Banach space $X$. Now setting $r=(1-2\epsilon)\sqrt{\log n}/\sqrt{n}$ and then going by the similar lines of argument as in \cite[pp. 743-744]{Bay}, it can be proved that the above inequality will be satisfied for $n$ large enough. As a result, $R_{p, q}^n(X)\geq(1-2\epsilon)\sqrt{\log n}/\sqrt{n}$ for sufficiently large $n$, and therefore 
$$
\liminf_{n\to\infty}R_{p, q}^n(X)\sqrt{n}/\sqrt{\log n}\geq 1.
$$
We skip the technical details of the aforesaid proof to avoid repetition.
Also, for the proof of the fact that $\limsup_{n\to\infty}R_{p, q}^n(X)\sqrt{n}/\sqrt{\log n}\leq 1$, take the complex-valued polynomial-say, $p(z)$-as given in \cite[p. 745]{Bay} (see also \cite{Boas}) and construct the $X$-valued polynomial $p_1(z)=p(z)e$ defined on $\ID^n$, where $e$ is a unit vector in $X$. Now we can adopt exactly the same lines of argument as in \cite[p. 745]{Bay} to validate our assertion.
\epf
\bpf[\bf Proof of Theorem \ref{N2thm4}]
(i) Given  a complex-valued holomorphic function $f$ with an expansion $(\ref{N2eq2})$ in $\ID^n$ (\lq$x_\alpha$\rq s are complex numbers in this case) and satisfying $\|f\|_{H^\infty(\ID^n, \IC)}\leq 1$, 
an application of H\"older's inequality yields
\begin{align*}
|x_0|^p+\sum_{k=1}^\infty r^{kp}\sum_{|\alpha|=k}|x_\alpha|^p
&=\sum_{k=0}^\infty \sum_{|\alpha|=k}|x_\alpha|^{2-p}r^{kp}|x_\alpha|^{2p-2}\\
&\leq\left(\sum_{k=0}^\infty r^{\frac{kp}{2-p}}\sum_{|\alpha|=k}|x_\alpha|\right)^{2-p}\left(\sum_{k=0}^\infty \sum_{|\alpha|=k}|x_\alpha|^2\right)^{p-1}\\
&\leq\left(\sum_{k=0}^\infty r^{\frac{kp}{2-p}}\sum_{|\alpha|=k}|x_\alpha|\right)^{2-p}.
\end{align*}
Therefore, $r_p^n(\IC)\geq (r_1^n(\IC))^{(2-p)/p}$.
Since $\lim_{n\to\infty}r_1^n(\IC)\left(\sqrt{n}/\sqrt{\log n}\right)=1$
(cf. \cite{Bay}), we have
$$
\liminf_{n\to\infty}
r_p^n(\IC)\left(\frac{n}{\log n}\right)^{\frac{2-p}{2p}}\geq 
\liminf_{n\to\infty}\left(r_1^n(\IC)\sqrt{\frac{n}{\log n}}\right)^{\frac{2-p}{p}}
=1,
$$
and thus $r_p^n(\IC)\geq C((\log n)/n)^{(2-p)/2p}$ for some constant $C>0$ and for all $n>1$.
The upper bound 
$r_p^n(\IC)\leq D\left((\log n)/n\right)^{(2-p)/2p}$ for some $D>0$ has already been established in \cite[p. 76]{Dja}. This completes the proof.

\vspace{3pt}
\noindent
(ii) To handle the second part of this theorem, we first construct $g(u)$ as in $(\ref{N2eq23})$ from a given holomorphic $f:\ID^n\to X$ with an expansion $(\ref{N2eq2})$ and satisfying $\|f\|_{H^\infty(\ID^n, X)}\leq 1$. Now, since $X$ is $p$-uniformly $PL$-convex, from the proof of \cite[Proposition 2.1(ii)]{Pav} we obtain
$$
\|P_1(z_0)\|\leq\frac{2}{(I_p(X))^{\frac{1}{p}}}(1-\|x_0\|^p)^{\frac{1}{p}}
$$
for any arbitrary $z_0\in\mathbb{T}^n$. Using a standard averaging trick (see f.i. \cite[p. 94]{Bla2}), it can be shown that the $P_1(z_0)$ in the above inequality could be replaced by $P_k(z_0)$ for any $k\geq 2$. Thus, we conclude
\be\label{N2eq20}
\sup_{z\in\mathbb{T}^n}\|P_k(z)\|
\leq\frac{2}{(I_p(X))^{\frac{1}{p}}}(1-\|x_0\|^p)^{\frac{1}{p}}.
\ee
Now, from \cite[Lemma 25.18]{Def3} it is known that there exists $R>0$ such that
$$
\left(\sum_{|\alpha|=k}\|x_\alpha\|^p\right)R^{kp}\leq\int_{\mathbb{T}^n}\|P_k(z)\|^pdz.
$$
Using inequality $(\ref{N2eq20})$ gives
$$
\sum_{|\alpha|=k}\|x_\alpha\|^p\leq\frac{2^p}{I_p(X)R^{kp}}(1-\|x_0\|^p).
$$
Assuming $r<R$, it is easy to see that
\begin{align*}
\|x_0\|^p+\sum_{k=1}^\infty r^{kp}\sum_{|\alpha|=k}\|x_\alpha\|^p
&\leq\|x_0\|^p+\frac{2^p}{I_p(X)}(1-\|x_0\|^p)\sum_{k=1}^\infty\left(\frac{r}{R}\right)^{kp}\\
&\leq\|x_0\|^p+\frac{2^p}{I_p(X)}(1-\|x_0\|^p)\frac{r^p}{R^p-r^p},
\end{align*}
which is less than or equal to $1$ if 
$$
r\leq R\left(\frac{I_p(X)}{2^p+I_p(X)}\right)^{\frac{1}{p}}=\left(\frac{I_p(X)}{2^p+I_p(X)}\right)^{\frac{2}{p}},
$$
as from the arguments in \cite[p. 627]{Def3}, it is clear that we can take
$R^p=I_p(X)/(I_p(X)+2^p)$. This proves the lower estimate for $r_p^n(X)$, and the upper estimate is trivial due to the fact that $r_p^n(X)\leq r_p^n(\IC)=1$ for $p\geq 2$.
\epf
\bpf[\bf Proof of Theorem \ref{N2thm5}]
Given any holomorphic $f:\ID^n\to\IC$ with an expansion $(\ref{N2eq2})$, we construct $g:\ID\to\IC$ from $f$ as in $(\ref{N2eq23})$. Provided $f(0)=1$ and $\mbox{Re}(f(z))>0$ for all $z\in\ID^n$, it is clear that $\mbox{Re}(g(u))>0$ for all $u\in\ID$ and $g(0)=1$. From Carath\'eodory's inequality, $|P_k(z)|\leq 2$ for all $z\in\mathbb{T}^n$, where $P_k(z)=\sum_{|\alpha|=k}x_\alpha z^\alpha$, $x_\alpha\in\IC$. Using this information, we will now give proofs for the cases $p\geq 2$ and $p<2$ separately.

\vspace{3pt}
\noindent
\underline{Case $p\in [2, \infty)$}: It is evident that
$$
\left(\sum_{|\alpha|=k}|x_\alpha|^p\right)^{\frac{1}{p}}\leq\left(\sum_{|\alpha|=k}|x_\alpha|^2\right)^{\frac{1}{2}}=\left(\int_{\mathbb{T}^n}|P_k(z)|^2dz\right)^{\frac{1}{2}}\leq 2
$$
for $p\geq 2$. A little calculation using this inequality shows that for $z\in r\ID^n$
$$
\frac{1}{2}\left(|x_0|^p+\sum_{k=1}^\infty\sum_{|\alpha|=k}|x_\alpha z^\alpha|^p\right)^{\frac{1}{p}}\leq
\frac{1}{2}\left(1+\frac{2^pr^p}{1-r^p}\right)^{\frac{1}{p}},
$$
which is less than or equal to $1$ if $r\leq \left((2^p-1)/(2^{p+1}-1)\right)^{1/p}$, i.e. 
$$
H_p^n\geq\left(\frac{2^p-1}{2^{p+1}-1}\right)^{\frac{1}{p}}.
$$
The reverse inequality follows from the fact that 
$$
H_p^n(\xi)=\left(\frac{2^p-1}{2^{p+1}-1}\right)^{\frac{1}{p}},
$$
where $\xi(z)=(1+z_1)/(1-z_1)$, $z=(z_1, z_2, \cdots, z_n)\in\ID^n$.

\noindent
\vspace{3pt}
\underline{Case $p\in(0,2)$}: We take care of the upper bound first. Indeed, for any holomorphic $f:\ID^n\to\ID$ with $f(0)=0$, the function $F:\ID^n\to\IC$ given by $F(z)=1-f(z)$ satisfies $\mbox{Re}(F(z))>0$ and $F(0)=1$. Therefore, we take $f$ as the function given in \cite[p. 76]{Dja} and follow the method illustrated there to show that 
$$
H_p^n\leq D\left(\frac{\log n}{n}\right)^{\frac{2-p}{2p}}
$$ 
for some constant $D>0$. The proof of the lower bound requires more work. Note that the hypercontractivity of the polynomial Bohnenblust-Hille inequality (see for example \cite[Theorem 8.19]{Def3}) ensures the existence of a constant $M>0$ such that 
\be\label{N2eq24}
\left(\sum_{|\alpha|=k}|x_\alpha|^{\frac{2k}{k+1}}\right)^{\frac{k+1}{2k}}\leq M^k\sup_{z\in\mathbb{T}^n}|P_k(z)|\leq 2M^k
\ee
for all $k\in\IN$. Using this inequality, we tackle the problem of finding a lower bound for $H_p^n$ by dividing it into two subcases.

\noindent
\underline{Subcase $p\in(0,1]$}: 
It is clear that we only need to find an $r>0$ such that the inequality 
\be\label{N2eq25}
\sum_{k=1}^\infty r^{kp}\sum_{|\alpha|=k}|x_\alpha|^p\leq 2^p-1
\ee
is satisfied.
Making use of the inequality $(\ref{N2eq24})$
and the fact that
$$
{n+k-1 \choose k}<e^k\left(1+\frac{n}{k}\right)^k,
$$
we get for any $r>0$:
\be\label{N2eq26}
\begin{aligned}
\sum_{k=1}^\infty r^{kp}\sum_{|\alpha|=k}|x_\alpha|^p &\leq\sum_{k=1}^\infty r^{kp}\left(\sum_{|\alpha|=k}|x_\alpha|^{\frac{2k}{k+1}}\right)^{\frac{p(k+1)}{2k}}{n+k-1 \choose k}^{1-\frac{p(k+1)}{2k}}\\
&\leq 2^p\sum_{k=1}^\infty M^{kp} r^{kp}{n+k-1 \choose k}^{1-\frac{p(k+1)}{2k}}\\
&\leq 2^p\sum_{k=1}^\infty e^k M^{kp} r^{kp}\left(\left(1+\frac{n}{k}\right)^{\frac{1}{p}-\frac{k+1}{2k}}\right)^{kp}.
\end{aligned}
\ee
Now observe that
$$
\left(1+\frac{n}{k}\right)^{\frac{1}{p}-\frac{k+1}{2k}}
\leq 2^{\frac{1}{p}-\frac{k+1}{2k}}\max\left\{1, \left(\frac{n}{k}\right)^{\frac{1}{p}-\frac{k+1}{2k}}\right\}\leq 2^{\frac{1}{p}-\frac{1}{2}}\max\left\{1, \left(\frac{n}{k}\right)^{\frac{1}{p}-\frac{k+1}{2k}}\right\},
$$
and that
$$
\left(\frac{n}{k}\right)^{\frac{1}{p}-\frac{k+1}{2k}}
=\left(\frac{n}{k}\right)^{\frac{1}{p}-\frac{1}{2}}\frac{k^{\frac{1}{2k}}}{n^{\frac{1}{2k}}}
\leq \left(\frac{n}{k}\right)^{\frac{1}{p}-\frac{1}{2}}\frac{e^{\frac{1}{2e}}}{n^{\frac{1}{2k}}},
$$
as the function $p(x):=x^{1/x}:(0, \infty)\to\IR$ attains its maximum at $x=e$. Further, by means of some calculations, we get that the following inequality is satisfied:
$$
q(x):=\frac{\left(\frac{n}{x}\right)^{\frac{1}{p}-\frac{1}{2}}}{n^\frac{1}{2x}}\leq q\left(\frac{p\log n}{2-p}\right)=M_1\left(\frac{n}{\log n}\right)^{\frac{2-p}{2p}}
$$
for some constant $M_1>0$, where $q:(0, \infty)\to\IR$ is differentiable. This asserts that 
$$
\left(1+\frac{n}{k}\right)^{\frac{1}{p}-\frac{k+1}{2k}}
\leq M_2\left(\frac{n}{\log n}\right)^{\frac{2-p}{2p}}
$$
for another constant $M_2>0$, and therefore, from $(\ref{N2eq26})$ we have
$$
\sum_{k=1}^\infty r^{kp}\sum_{|\alpha|=k}|x_\alpha|^p\leq 2^p\sum_{k=1}^\infty\left(M_3r\left(\frac{n}{\log n}\right)^{\frac{2-p}{2p}}\right)^{kp}
$$
for a new constant $M_3>0$. Clearly, the inequality $(\ref{N2eq25})$ is satisfied if we take $r=\gamma\left((\log n)/n\right)^{(2-p)/2p}$ for some sufficiently small constant $\gamma>0$. This finishes the proof for this part.

\vspace{3pt}
\noindent
\underline{Subcase $p\in(1, 2)$}: Proof for this part is in fact similar to the proof of Theorem \ref{N2thm4}(i) and also to the proof of \cite[Proposition 1.4]{Bla2}. For any $z\in r\ID^n$ where 
$$
r=(H_1^n)^{\frac{2-p}{p}}(H_2^n)^{\frac{2p-2}{p}},
$$ 
we have
\begin{align*}
|x_0|^p+\sum_{k=1}^\infty\sum_{|\alpha|=k}|x_\alpha z^\alpha|^p &\leq\sum_{k=0}^\infty\sum_{|\alpha|=k}|x_\alpha|^{2-p}(H_1^n)^{k(2-p)}|x_\alpha|^{2p-2}(H_2^n)^{k(2p-2)}\\
&\leq\left(\sum_{k=0}^\infty\sum_{|\alpha|=k}|x_\alpha|(H_1^n)^{k}\right)^{2-p}\left(\sum_{k=0}^\infty\sum_{|\alpha|=k}|x_\alpha|^{2}(H_2^n)^{2k}\right)^{p-1}\\
&\leq (2^{2-p})(2^{2(p-1)})=2^p.
\end{align*}
Consequently, using the already proven estimates for $H_1^n$ and $H_2^n$, we obtain
$$
H_p^n\geq (H_1^n)^{\frac{2-p}{p}}(H_2^n)^{\frac{2p-2}{p}}\geq \left(C_1\sqrt{\frac{\log n}{n}}\right)^{\frac{2-p}{p}}\left(\frac{3}{7}\right)^{\frac{p-1}{p}}= C\left(\frac{\log n}{n}\right)^{\frac{2-p}{2p}}
$$
for some constant $C>0$. Hence the proof is complete.
\epf
\section{An additional observation}\label{N2sec3}
Let $G$ be a domain in the complex Banach space $Y$. For
any holomorphic mapping $f:G\to X$, let $D^kf(y)$ denote the $k$-th Fr\'echet derivative $(k\in\IN)$
of $f$ at $y\in G$, which is a bounded symmetric $k$-linear mapping
from $\prod_{i=1}^k Y\to X$. Any such $f$ can be expanded into the series
\be\label{N2eq21}
f(y)=\sum_{k=0}^\infty\frac{1}{k!}D^kf(y_0)((y-y_0)^k)
\ee
in a neighborhood of any given $y_0\in G$. It is understood that
$D^0f(y_0)(y^0)=f(y_0)$ and
$$
D^kf(y_0)(y^k)=D^kf(y_0)(\underbrace{y,y,\cdots, y}_\text{$k$-times})
$$
for $k\geq 1$. Now for $1\leq p, q<\infty$, we define
$$
R_{p, q}(f, Y, X)=\sup\left\{r\geq 0:\|x_0\|^p+\left(\sum_{k=1}^\infty\frac{r^k}{k!}\left\|D^kf(0)\right\|\right)^q\leq 1\right\}
$$
for any holomorphic $f$ defined in the open unit ball $B$ of $Y$ with an expansion $(\ref{N2eq21})$ in a neighborhood of $y_0=0$, while $f(0)=x_0$ and $\|f(y)\|\leq 1$ for all $y\in B$. 
We mention here that for any $k$-linear mapping $L:\prod_{i=1}^k Y\to X$,
$$
\|L\|:=\sup_{\substack{\|y_i\|\leq 1\\ i=1, 2\cdots k}}\|L(y_1, y_2, \cdots, y_k)\|.
$$
We further define 
$$
R_{p, q}(Y, X)=\{R_{p, q}(f, Y, X): f:B\to X \,\mbox{holomorphic and}\, \|f(y)\|\leq 1 \,\mbox{for all}\, y\in B\}.
$$
After all these preparations, we will now show that the existence of nonzero Bohr radius $R_{p, q}(Y, X)$ depends only on $X$.
\bprop\label{N2prop2}
Given any two complex Banach spaces $X$ and $Y$, $R_{p, q}(Y, X)>0$ if and only if $R_{p, q}(X)>0$.
\eprop
\bpf
Suppose $R_{p, q}(X)>0$. Then for any fixed $y_1\in B$, we construct the holomorphic function $g_1:\ID\to X$ with $\|g_1(u)\|\leq 1$ for all $u\in\ID$; given by
$$
g_1(u)=f(uy_1)=x_0+\sum_{k=1}^\infty\left(\frac{1}{k!}D^kf(0)(y_1^k)\right)u^k.
$$
It is clear that for any $0<|u|=R\leq R_{p, q}(X)$ we have
$$
\|x_0\|^p+\left(\sum_{k=1}^\infty \frac{R^k}{k!}\left\|D^kf(0)(y_1^k)\right\|\right)^q\leq 1,
$$
and since $y_1$ is arbitrary, this gives (see f.i. \cite[Remark 15.19]{Def3}):
\begin{align*}
\frac{1}{k!}\left\|D^kf(0)\right\|\leq\frac{k^k}{(k!)^2}\sup_{y\in B}\left\|D^kf(0)(y^k)\right\|&= \frac{k^k}{k!}\sup_{y\in B}\frac{1}{k!}\left\|D^kf(0)(y^k)\right\|\\
&\leq\frac{k^k}{k!R^k}(1-\|x_0\|^p)^{\frac{1}{q}}\leq\frac{e^k}{R^k}(1-\|x_0\|^p)^{\frac{1}{q}}.
\end{align*}
Using this inequality we have, after a little calculation, that $R_{p, q}(Y, X)\geq R/(2e)>0$. Conversely, assume now $R_{p, q}(Y, X)>0$. For a given $\delta>0$, we take any $x_0\in X$ such that $\|x+zx_0\|\leq 1+\delta$
for all $z\in\ID$ and for all $x\in X$ with $\|x\|=1$. We now define 
the holomorphic function $F:B\subset Y\to X$ by $F(y)=x+l(y)x_0$ where $l\in X^*$ with $\|l\|=1$. Thus $\|F(y)\|\leq 1+\delta$ for all $y\in B$. Noting that $DF(0)(y)=l(y)x_0$ and that $\|DF(0)\|=\|x_0\|$, we have
$$
\left(\frac{1}{1+\delta}\right)^p+\frac{R^q\|x_0\|^q}{(1+\delta)^q}\leq 1
$$
for any $R\leq R_{p, q}(Y, X)$. Consequently, $\|x_0\|\leq (1/R)\left((1+\delta)^q-(1+\delta)^{q-p}\right)^{1/q}$,
which then implies inequality $(\ref{N2eq6})$. As a consequence, $R_{p, q}(X)>0$.
\epf
\section*{Acknowledgements}
The author is thankful to Prof. Bappaditya Bhowmik for his kind help in obtaining a softcopy of \cite{Dja}. He also thanks 
Mr. Aritra Bhowmick for some stimulating discussions during the preparation of this manuscript.

\end{document}